\documentclass[twocolumn]{amsart}
\usepackage{amsrefs,tikz,caption,times}
\usetikzlibrary{decorations}
\usepgflibrary{snakes}
\newtheorem*{theorem}{Theorem}
\begin{document}

\title{Orange peels and Fresnel Integrals}
\date{February 6, 2012}
\author{Laurent Bartholdi}
\address{Mathematisches Institut, Georg-August Universit\"at zu G\"ottingen}
\email{laurent@uni-math.gwdg.de}

\author{Andr\'e G. Henriques}
\address{Mathematisch Instituut, Universiteit Utrecht}
\email{A.G.henriques@uu.nl}

\thanks{Partially supported by the Courant Research Centre ``Higher
  Order Structures'' of the University of G\"ottingen}

\maketitle

{\bf
There are two standard ways of peeling an orange: either cut the skin
along meridians, or cut it along a spiral.  We consider here the
second method, and study the shape of the spiral strip, when unfolded
on a table.  We derive a formula that describes the corresponding
flattened-out spiral.  Cutting the peel with progressively thinner
strip widths, we obtain a sequence of increasingly long spirals.  We
show that, after rescaling, these spirals tends to a definite shape,
known as the Euler spiral.  The Euler spiral has applications in many
fields of science.  In optics, the illumination intensity at a point
behind a slit is computed from the distance between two points on the
Euler spiral.  The Euler spiral also provides optimal curvature for
train tracks between a straight run and an upcoming bend.  It is
striking that it can be also obtained with an orange and a kitchen
knife.}

\section{Outline}
\begin{figure}[h]
\[
\begin{tikzpicture}
  \node[scale = .925] at (-.4,0) {\includegraphics[scale=0.15]{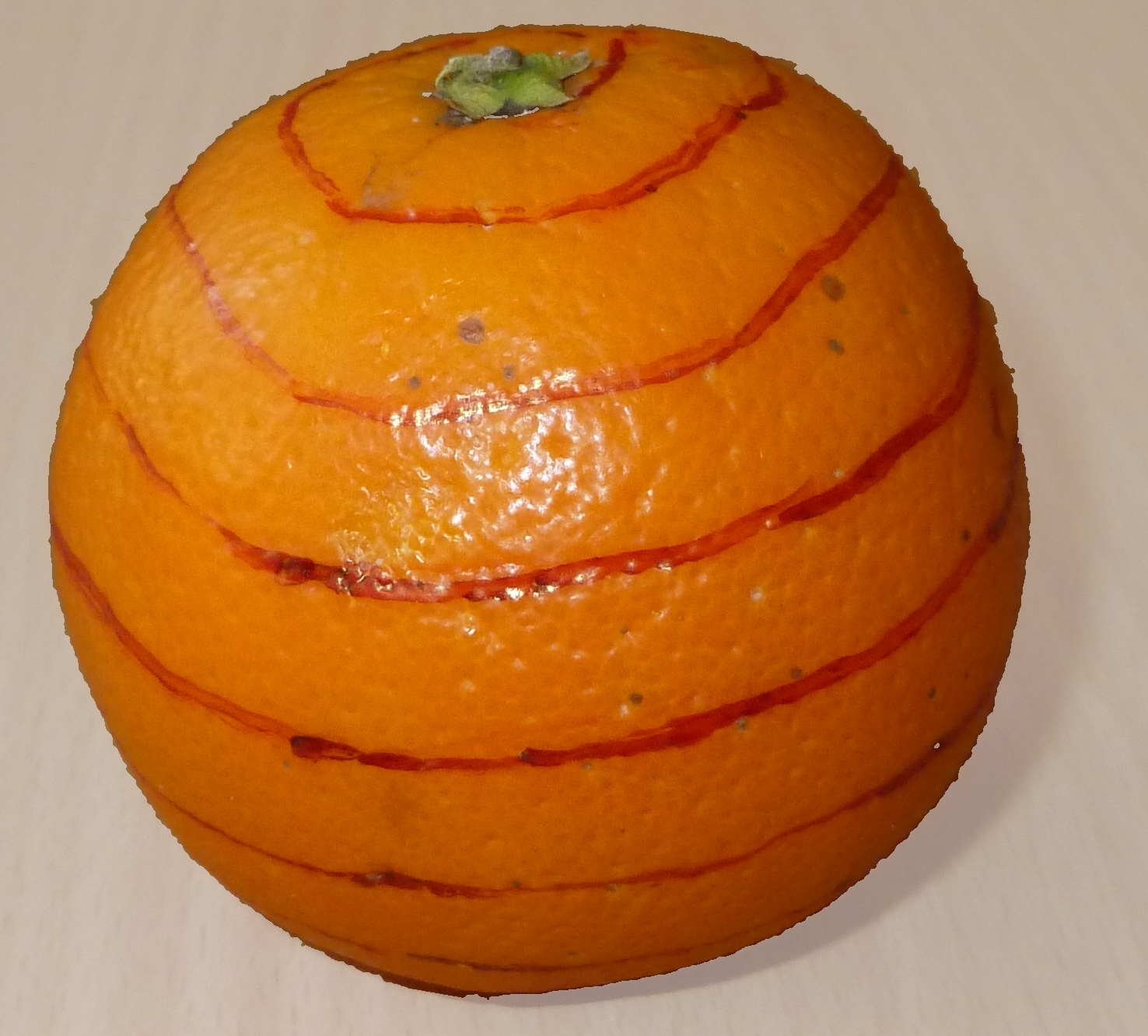}};
  \draw[<->,white,very thick] (0.95,1.45) -- node[right, pos=.3, scale=.97] {$1/N$} (1.6,0.4);
\end{tikzpicture}
\]
\caption{An orange, assumed to be a sphere of radius one, with spiral of width $1/N$.\vspace{-1cm}}\label{fig:orange}
\end{figure}

\begin{figure}[h]
\[
\begin{tikzpicture}
  \node at (0,0) {\includegraphics[scale=0.098]{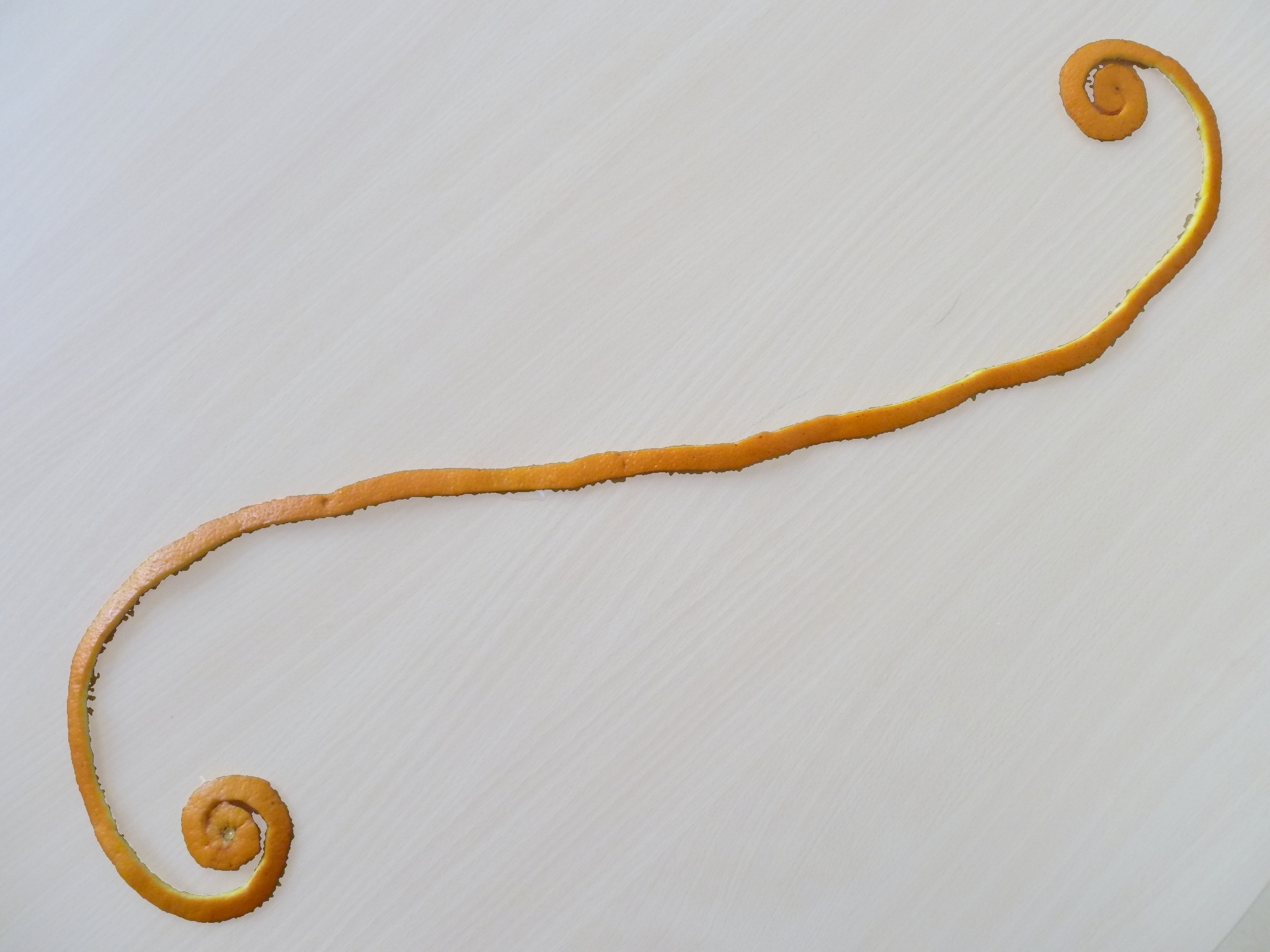}};
  \draw[very thick] (2.5,-2) -- (2.5,-2.05) -- node[below, scale=.97] {1\,cm} (2.8,-2.05) -- (2.8,-2);
\end{tikzpicture}
\]
\caption{The flattened-out orange peel.}\label{fig:peel}
\end{figure}

Cut the skin of an orange along a thin spiral of constant width
(fig.~\ref{fig:orange}) and lay it flat on a table
(fig.~\ref{fig:peel}).  A natural breakfast question, for a
mathematician, is which shape the spiral peel will have, when
flattened out.  We derive a formula that, for a given cut width,
describes the corresponding spiral's shape.

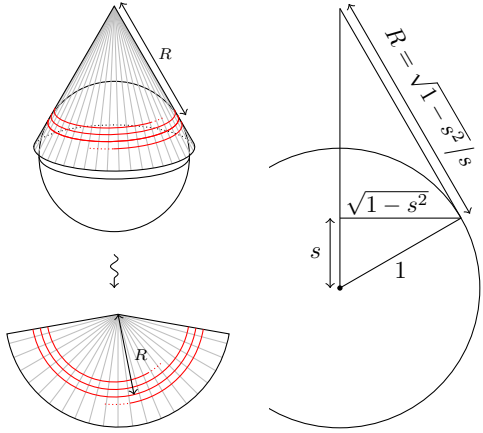
\begin{figure}[h]
\[
  \begin{tikzpicture}
  \draw (0,0) circle (1);
  \coordinate (c) at (0,2);
 \begin{scope}
  \pgftransformcm{1}{0}{0}{.3}{\pgfpoint{0}{0}}
  \draw (-1,0) arc (-180:0:1); 
  \pgftransformcm{1.07}{0}{0}{1.07}{\pgfpoint{0}{12}}
  \draw[densely dotted] (-1,0) arc (175:5:1); 
  \foreach \x in {180, 188.5, ..., 360}{\draw[gray!50] (\x:1) -- (c);}
  \draw (-190:1) arc (-190:10:1) -- (c) -- cycle; 
  \pgftransformcm{.9}{0}{0}{.9}{\pgfpoint{0}{16}}
  \pgftransformcm{.95}{.1}{0}{1}{\pgfpoint{-.4}{8.2}}
  \draw[red, densely dotted] (-110:1) arc (-110:-90:1); 
  \draw[red] (-90:1) arc (-90:10:1); 
  \pgftransformcm{.95}{0}{0}{.95}{\pgfpoint{0}{8.2}}
  \draw[red] (-190:1) arc (-190:10:1); 
  \draw[red] (-190:1) arc (-190:10:1); 
  \pgftransformcm{.95}{0}{0}{.95}{\pgfpoint{0}{8.2}}
  \draw[red] (-190:1) arc (-190:10:1); 
  \pgftransformcm{.95}{0}{0}{.95}{\pgfpoint{0}{8.2}}
  \draw[red] (-190:1) arc (-190:-50:1); 
  \draw[red, densely dotted] (-50:1) arc (-50:-30:1); 
 \end{scope}
  \draw[<->] (0,2)+(30:.1) -- node[scale=.85, right, pos=.45]{$\scriptstyle R$} (30:1.1);
  \draw[->,decorate,decoration={snake,post length=.5mm,amplitude=.5mm,segment length=2mm}] (0,-1.3) -- (0,-1.75);
  \begin{scope}[xshift=-3.95cm,yshift=-2.6cm]
  \foreach \x in {-15, -22.5, ..., -165}{\draw [gray!50] (4,.5) -- +(\x:1.5);}
  \draw (4,.5) -- +(-10:1.5) arc (-10:-170:1.5) -- cycle;
  \draw [red] (3.95,.5) +(-80:1.2) arc (-80:-9.5:1.2);
  \draw [red] (3.95,.5) +(-169.5:1.1) arc (-169.5:-9.5:1.1);
  \draw [red] (3.95,.5) +(-169.5:1.0) arc (-169.5:-9.5:1.0);
  \draw [red] (3.95,.5) +(-169.5:.9) arc (-169.5:-60:.9);
  \draw[red, densely dotted] (3.95,.5) +(-60:.9) arc (-60:-45:.9);
  \draw[red, densely dotted] (3.95,.5) +(-96:1.2) arc (-96:-80:1.2);
  \draw[<->] (4,.5) -- node[scale=.85, right, pos=.5]{$\scriptstyle R$} +(-78.7:1.09);
  \end{scope}
  \begin{scope}[xshift=3cm,yshift=-1.75cm]
  \pgftransformcm{.93}{0}{0}{.93}{\pgfpoint{0}{0}}
  \clip (-1,-2.02) rectangle (2.1,4.2);
  \draw (0,0) circle (2cm);
  \draw (0,0) -- (0,1);
  \draw[<->] (-.14,0) -- node[left, scale=.93] {$s$} (-.14,1);
  \draw (0,1) -- (0,4);
  \draw (0,1) -- node[above, xshift=-5,yshift=-2, scale=.9] {$\sqrt{1-s^2}$} (30:2);
  \draw (0,0) -- node[below, scale=.93] {$1$} (30:2);
  \draw (30:2) -- (0,4);  
  \filldraw (0,0) circle (.03); 
  \pgftransformxshift{3.6}
  \pgftransformyshift{1.8}
  \draw[<->] (30:2) -- node[above,sloped, scale=.9] {$R=\sqrt{1-s^2\,}\!\big/s$} (0,4);
  \end{scope}
\end{tikzpicture}
\]
\caption{\emph{left:} Spiral on the sphere, transferred to the tangent cone, and developed
on the plane, for computing its radius of curvature;\\
  \emph{right:} The computation of the radius of curvature $R$ of the flattened spiral.}
\label{fig:cone}
\end{figure}

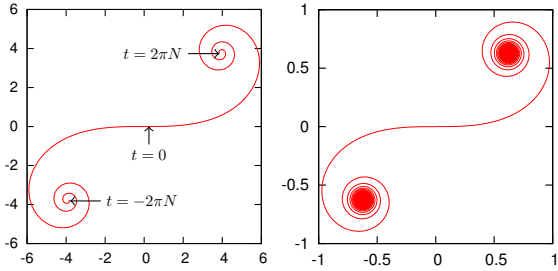
\begin{figure}[h]
\[
\begin{tikzpicture}[scale=0.8]
 \useasboundingbox (-2.2,-2.1) rectangle (1.7,1.9);
 \node at (-0.3,-0.1) {\includegraphics[scale=0.8]{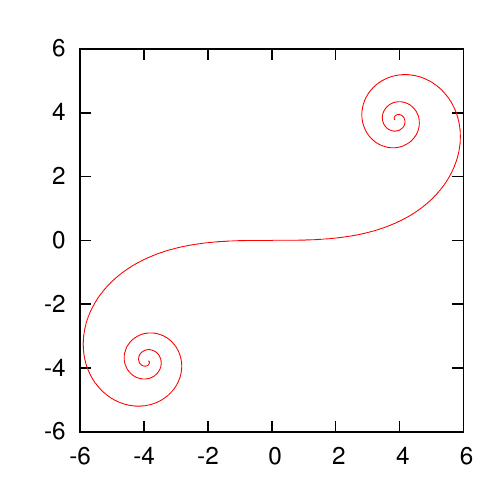}};
 \draw[->] (0,-.3) node[anchor=north, scale=.6] {$t=0$} -- (0,0);
 \draw[->] (-0.8,-1.24) node[anchor=west, scale=.6] {$t=-2\pi N$} --
(-1.31,-1.24);
 \draw[->] (0.65,1.21) node[anchor=east, scale=.6] {$t=2\pi N$} -- (1.17,1.21);
\end{tikzpicture}\quad\,\,\,
\begin{tikzpicture}[scale=0.865]
 \useasboundingbox (-2.3,-2.0) rectangle (1.6,1.8);
 \pgftransformyshift {-1.7}
 \node at (-0.4445,-0.1) {\includegraphics[scale=0.865]{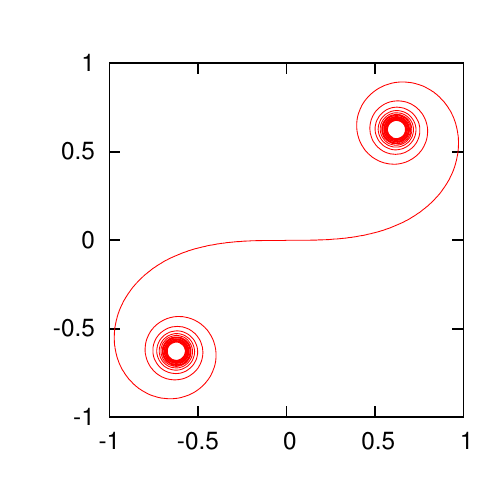}};
 \fill[red] (1.05,1.12) circle (.1) (-1.19,-1.13) circle (.1);
\end{tikzpicture}
\]
\caption{\emph{left:} Maple plot of the orange peel spiral ($N=3$);\\
\emph{right:} the Euler spiral; limit $N \to \infty$.}\label{fig:spiral}
\end{figure}

For the analysis, we parametrize the spiral curve by a constant speed
trajectory, and express the curvature of the flattened-out spiral as a
function of time. 
This is achieved by comparing a revolution of the
spiral on the orange with a corresponding spiral on a cone tangent to
the surface of the orange (fig.~\ref{fig:cone}, left).  Once we know the
curvature, we derive a differential equation for our spiral, which we
solve analytically (fig.~\ref{fig:spiral}, left).

We then consider what happens to our spirals when we vary the strip
width.  Two properties are affected: the overall size, and the shape.
Taking finer and finer widths of strip, we obtain a sequence of
increasingly long spirals; rescale these spirals to make them all of
the same size.  We show that, after rescaling, the shape of these
spirals tends to a well defined limit.  The limit shape is a classical
mathematical curve, known as the \emph{Euler spiral} or the
\emph{Cornu spiral} (fig.~\ref{fig:spiral}, right).  This spiral is the
solution of the \emph{Fresnel integrals}.

The Euler spiral has many applications.  In optics, it occurs in the
study of light diffracting through a
slit~\cite{hecht:optics}*{\S10.3.8}.  More precisely, the illumination
intensity at a point behind a slit is the square of the distance
between two points on the Euler spiral, easily determined from the
slit's geometry.

The same spiral is also used in civil engineering: it provides optimal
curvature for train tracks~\cite{profillidis:railway}*{\S14.1.2}.  A
train that travels at constant speed and increases the curvature of
its trajectory at a constant rate will naturally follow an arc of the
Euler spiral.  The review~\cite{Levien:EECS-2008-111} describes the
history of the Euler spiral and its three independent discoveries.

\section{Analysis}
For the purpose of our mathematical treatment, we shall replace the
orange by a sphere of radius one.  The spiral on the sphere is taken
of width $1/N$, see (fig.~\ref{fig:orange}).  The area of the sphere
is $4\pi$, so the spiral has a length of roughly $4\pi N$.  We
describe the flattened-out orange peel spiral by a curve $(x(t),y(t))$
in the plane, parameterized at unit-speed from time $t=-2\pi N$ to
$t=2\pi N$.

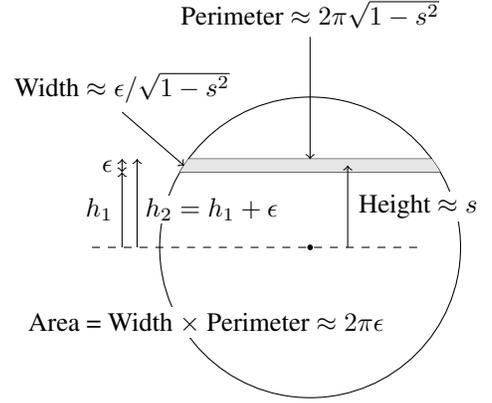
\begin{figure}[h]
\[
\begin{tikzpicture}
  \draw (0,0) circle (2cm);
  \draw (30:2) -- (150:2);
  \draw (36:2) -- (144:2);
  \fill[gray!20!white] (150:2) arc (150:144:2cm) -- (36:2) arc (36:30:2cm) -- cycle;
  \draw[->] (0,2.8) node[anchor=south] {Perimeter $\approx 2\pi\sqrt{1-s^2}$} -- (0,0 |- -36:-2);
  \draw[->] (-2.5,1.8) node[anchor=south] {Width $\approx \epsilon/\sqrt{1-s^2}$} -- (147:2);
  \draw[->] (-2.5,0) -- node[left] {$h_1$} (-2.5,0 |- 30:2);
  \draw[<->] (-2.5,0 |- 30:2) -- node[left] {$\epsilon$} (-2.5,0 |- 36:2);
  \draw[->] (-2.3,0) -- node[right, yshift=-2.5, xshift=2, fill=white, inner sep=1] {$h_2=h_1+\epsilon$} (-2.3,0 |- 36:2);
  \draw[dashed] (-2.9,0) -- (1.5,0);
  \draw[->] (.5,0) --node[anchor=west, fill=white, inner sep=1, xshift=3] {Height $\approx s$} (.5,0 |- 33:2);
  \node[fill=white, inner sep=2] at (-1.38,-1) {Area = Width $\times$ Perimeter $\approx 2\pi\epsilon$};
  \filldraw (0,0) circle (.03); 
\end{tikzpicture}
\]
\caption{Area of a thin circular strip on the sphere.}
\label{fig:area}
\end{figure}

On a sphere of radius one, the area between two horizontal planes at
heights $h_1$ and $h_2$ is ${2\pi(h_2-h_1)}$, see (fig.~\ref{fig:area}).  It
follows that, at time $t$, the point on the sphere has height
$s:=t/2\pi N$.

Our first goal is to find a differential equation for
$(x(t),y(t))$. For that, we compute the radius of curvature $R(t)$ of
the flattened-out spiral at time $t$: this is the radius of circle
with best contact to the curve at time $t$.  For example, $R(-2\pi
N)=R(2\pi N)=0$ at the poles, and $R(0)=\infty$ at the equator.

For $N$ large, the spiral at time $t$ follows roughly a parallel at
height $s$ on the orange. The surface of the sphere can be
approximated by a tangent cone whose development on the plane is a
disk sector (fig.~\ref{fig:cone}, left). The radius
\[
R(t)=\sqrt{1-s^2}/s=\sqrt{(2\pi N)^2-t^2}/t
\]
of that disk equals the radius of curvature of the spiral at time $t$,
and can be computed using Thales' theorem (fig.~\ref{fig:cone},
right). The radius $R(t)$ is in fact only determined up to sign; our
choice reflects the NE-SW orientation of the spiral on the sphere.


Now, the condition that we move move at unit speed on the sphere ---
and on the plane --- is $(\dot x)^2+(\dot y)^2=1$, and the condition that
the spiral has a curvature of $R(t)$ is $\dot x\ddot y-\ddot x\dot y=1/R$.
Here, $\dot x$ and $\dot y$ are the speeds of $x$ and $y$ respectively, and $\ddot x$ and $\ddot y$ are their accelerations.
In fact, introducing the complex path $z(t)=x(t)+iy(t)$, the
conditions can be expressed as $|\dot z|^2=1$ and $\ddot z \dot {\bar z}=i/R$.


The solution has the general form \[z(t)=\int_0^t\exp(i\phi(u))du,\]
for a real function $\phi$; indeed, its derivative is computed as $\dot
z=\exp(i\phi(t))$ and has norm $1$. As $\ddot z\dot{\bar z}=i\dot\phi(t)$,
we have $\dot\phi(t)=s/\sqrt{1-s^2}$, which has as elementary solution
$\phi(t)=-\sqrt{(2\pi N)^2-t^2}$. We have deduced that the
flattened-out spiral has parameterization
\[\left\{\begin{array}{l}
    \displaystyle x(t)\,=\,\,\int_0^t\cos\sqrt{(2\pi N)^2-u^2}du,\\
    \displaystyle y(t)\,=\,-\int_0^t\sin\sqrt{(2\pi N)^2-u^2}du.
  \end{array}\right.\]
The flattened-out peel of an orange is shown in (fig.~\ref{fig:peel}),
and the corresponding analytic solution, computed by
\textsc{Maple}~\cite{Maple10}, is shown in (fig.~\ref{fig:spiral},
left). The orange's radius was 3cm, and the peel was 1cm wide, giving
$N=3$.

\section{Limiting behaviour}
What happens if $N$ tends to infinity, that is, if we peel the
orange with an ever thinner spiral?  For that, we recall the power
series approximation
\[
\sqrt{a^2-u^2}=a-\frac{u^2}{2a}+\mathcal O(\frac{u^4}{a^3}),
\]
which we substitute with $a=2\pi N$ in the above expression:
\begin{align*}
z(t)&=\int_0^t\exp\Big(\!-i\sqrt{(2\pi N)^2-u^2}\,\Big)du\\&
\approx\int_0^t\exp \Big(\!-i\Big(2\pi N-\frac{u^2}{2\cdot2\pi N}\Big)\!\Big)du.
\end{align*}
Taking only values of $N$ that are integers, this simplifies to
$\int_0^t\exp(iu^2/4\pi N)du$. We then set $v=u/\sqrt{4\pi N}$ to
obtain
\[
z(t)\approx \sqrt{4\pi N}\int_0^{t/\sqrt{4\pi N}}\exp(iv^2)dv.
\]
The approximation error is $\int_0^t\mathcal
O(\frac{u^4}{a^3})du=\mathcal O(t^5/N^3)$, which becomes negligible
compared to the size $\mathcal O(\sqrt N)$ of the spiral for $|t|\ll
N^{0.7}$.
   
The above curve is, up to scaling and parameterization speed,
the solution of the classical Fresnel integral
\[(X(t),Y(t))=\left(\int_0^t\cos u^2du,\int_0^t\sin u^2du\right),\]
defined by the condition that the radius of curvature at time $t$ is
$1/2t$; here the parameterization is over $t$ from $-\infty$ to
$+\infty$. The corresponding curve is called the Euler spiral and
winds infinitely often around the points
$\pm(\sqrt{\frac\pi8},\sqrt{\frac\pi8})$.  Setting $T:=t/\sqrt{4\pi
  N}$, the condition $|t|\ll N^{0.7}$ becomes $|T|\ll N^{0.2}$.  We
have thus proven:
\begin{theorem}
  If $T\ll N^{0.2}$, then the part of the orange peel of width $1/N$
  parameterized between $-\sqrt{4\pi N}\,T$ and $\sqrt{4\pi N}\,T$ is
  a good approximation for the part of the Euler spiral parameterized
  between $-T$ and $T$.
\end{theorem}

\noindent Note that for large $N$, the piece of the orange peel
parameterized between $-\sqrt{4\pi N}\,T$ and $\sqrt{4\pi N}\,T$ forms
a rather thin band around the orange's equator.  The contribution of
rest of the orange disappears due to the rescaling process.

\section{Conclusion}
The Euler spiral is a well known mathematical curve.  In this article,
we explained how to construct it with an orange and a kitchen knife.
Flattened fruit peels have already been considered, e.g.\ those of
apples~\cite{Turrell}, but were never studied analytically.  The Euler
spiral that we obtained has had many discoveries across
history~\cite{Levien:EECS-2008-111}; ours occurred over breakfast.

\end{document}